\documentclass[10pt,leqno]{amsart}

\usepackage[latin1]{inputenc}

\usepackage[a4paper,asymmetric]{geometry}
\usepackage{latexsym, amsmath, amssymb, esint}
\usepackage{bm}
\usepackage{graphics,epsfig,graphicx,graphics}
\usepackage{color, xcolor}
\usepackage{tikz}
\usepackage[hyperpageref]{backref}
\usetikzlibrary{patterns}
\usepackage[english]{babel}
\usepackage[colorlinks=true,urlcolor=blue, citecolor=red,linkcolor=blue,
linktocpage,pdfpagelabels, bookmarksnumbered,bookmarksopen]{hyperref}
\usepackage{hyperref}

\newtheorem{theorem}{Theorem}[section]

\newtheorem{lemma}[theorem]{Lemma}
\newtheorem{prop}[theorem]{Proposition}

\theoremstyle{definition}
\newtheorem{rem}[theorem]{Remark}
\newcommand{\pa}{\partial}

\newcommand{\deficit}{{\rm def}(\Omega)}

\numberwithin{equation}{section}

\begin{document}
\title[The method of moving planes: a quantitative approach]{The method of moving planes: a quantitative approach}

\author[G. Ciraolo and A. Roncoroni]{Giulio Ciraolo and Alberto Roncoroni}

\address{G. Ciraolo,  Dipartimento di Matematica e Informatica, Universit\`a degli Studi di Palermo, Via Archirafi 34, 90123, Palermo - Italy}
\email{giulio.ciraolo@unipa.it}

\address{A. Roncoroni, Dipartimento di Matematica F. Casorati, Universit\`a degli Studi di Pavia, Via Ferrata 5, 27100, Pavia - Italy}
\email{alberto.roncoroni01@universitadipavia.it}

\maketitle

\begin{abstract}
We review classical results where the method of the moving planes has been used to prove symmetry properties for overdetermined PDE's boundary value problems (such as Serrin's overdetermined problem) and for rigidity problems in geometric analysis (like Alexandrov soap bubble Theorem), and we give an overview of some recent results related to quantitative studies of the method of moving planes, where quantitative approximate symmetry results are obtained.

\end{abstract}

%
%	\begin{center}
%		\begin{minipage}{9cm}
%			\small
%			\tableofcontents
%		\end{minipage}
%	\end{center}
%

\bigskip



%%%%%%%%%%%%%%%%%%%%%%%%%%%%%%%%%%%%%%%%%%%%%%%%%
%%%%%%%%%%%%%%%%%%%%%%%%%%%%%%%%%%%%%%%%%%%%%%%%%

\section{
	Introduction}
The method of moving planes (MMP) was introduced by Alexandrov in \cite{Alex2} to prove what is nowadays called \emph{Alexandrov's soap bubble Theorem}, which asserts that spheres are the only connected closed embedded hypersurfaces in a space form (i.e. the Euclidean space, the Hyperbolic space and the hemisphere). Some years later, Serrin \cite{Serrin} employed the same method to prove a symmetry result in potential theory, which gave rise to the research field which is nowadays called \emph{overdetermined problems} for partial differential equations. 

Both Alexandrov's and Serrin's results originate a great interest in geometric analysis and PDE's communities. The MMP is a powerful tool which has been used to prove several results in geometric analysis, for elliptic and parabolic PDEs, Harnack's inequalities and many others (see e.g. \cite{Aronson_Caffarelli,Bennet1,Bennet2,Bonforte3,Bonforte1,Bonforte2,KKS,Li,Meeks,Reichel,Reichel1,Schoen}).
%(see e.g. \cite{KKS, Meeks, Schoen}\cite{Reichel1,Reichel,Aronson_Caffarelli,Bonforte1,Bonforte2,Bonforte3}\cite{Li}\cite{Bennet1,Bennet2})
One of the most influencing application in the theory of PDEs is the approach of Gidas Ni Nirenberg \cite{GNN1,GNN2} (see also \cite{BN1,BN2,BN3,Sciunzi2,Dancer,DP,DPR,DS,Li1,Li2,Monticelli,Sciunzi1}). The MMP was also used in \cite{CGS} to prove asymptotic radial symmetry of positive solutions for the conformal scalar curvature
equation and other semilinear elliptic equations (see also \cite{CL,Gidas,KMPS,Obata}).

The goal of this paper is to review the classical approach of the MMP (as Alexandrov and Serrin used) and to describe some recent results which are based on a quantitative version of the MMP.

Both the results of Alexandrov and Serrin apply the MMP and use maximum principles to obtain a symmetry result. The maximum principle is, in this setting, a tool to obtain qualitative information on a suitable function. In particular it is used to obtain a symmetry or rigidity result: the solution to a certain problem exists if and only if the domain satisfies some symmetry. Once the symmetry result is obtained, it is of interest to study its quantitative counterpart and to understand, in a quantitative way, whether or not the solution is close to the symmetry configuration. Hence, maximum principles must be replaced by quantitative tools (like for instance Harnack's type inequalities) and every step in the symmetry proof must be carefully quantified.

Now, we describe in more details the problems that we are going to consider. In order to simplify the exposition, we prefer to state the main theorems in the simplest possible setting, and add some remark regarding the generalizations available in literature. 

We start by the well-known Alexandrov soap bubble Theorem in geometric analysis on hypersurfaces with constant mean curvature $H$. We recall that the mean curvature $H$ of a hypersurface $S$ of class $C^2$ at $p\in S$ is given by the arithmetic mean of its principal curvatures at $p$.

\begin{theorem}[\cite{Alex_hyp,Alex2}] \label{thm_Alex}
	Let $\Omega \subset \mathbb{R}^n$ be a bounded connected domain with boundary $S=\partial \Omega$ of class $C^2$. Then the mean curvature $H$ of $S$ is constant if and only if $S$ is a sphere.
\end{theorem}

As already mentioned, Theorem \ref{thm_Alex} was proved by Alexandrov in \cite{Alex2} and it is the first paper where the method of moving planes appears (Alexandrov called it the Reflection Principle, the name method of moving planes appears in the 70s). We mention that since $S=\partial \Omega$ then $S$ is an embedded hypersurface; if this assumption is dropped then Theorem \ref{thm_Alex} does not hold in general (see \cite{HTY} and \cite{Wente} for classical counterexamples). For immersed hypersurfaces, an Alexandrov's type Theorem holds under the assumption that the hypersurface $S$ is of genus  $0$, i.e. $S$ is an immersed topological sphere (see \cite{Hopf} for the result and also \cite{Yau} for generalizations to higher dimensional hypersurfaces immersed in space forms and \cite{Nelli_survey} for more details) or by assuming that $S$ is stable in $\mathbb{R}^n$ (see \cite{BarbosaDoCarmo} and \cite{BdCE}). Moreover there exists non-closed constant mean curvature hypersufaces embedded in $\mathbb{R}^3$ which are not diffeomorphic to a sphere, like for instance the unduloids (see \cite{Delunay} and \cite{Wu_cinese} for the generalization to $\mathbb{R}^n$).

Still by using the MMP, several extensions of Theorem \ref{thm_Alex} have been proved in literature. For instance, it holds for hypersurfaces embedded in the Hyperbolic space and in the hemisphere (see \cite{Alex_hyp,Alex3} see also \cite{Hsiang-Hsiang}) and for other type of curvature functions, such as higher order curvature functions (see \cite{Alex3,CY,hartman,korevar,rosenberg}). Recently, and without using the MMP, new generalization of Alexandrov's Theorem \ref{thm_Alex} are given in \cite{Brendle}, where constant mean curvature hypersurfaces are studied in rotationally symmetric Riemannian manifolds (e.g. the space forms, the Schwarzschild, the DeSitter-Schwarzschild and Anti-DeSitter-Schwarzschild manifolds), and in \cite{DelgadinoMaggi} where the regularity assumptions on $\partial \Omega$ are minimal. 

%\todo{anche \cite{AIR,AD} - non le metterei...brendle se le mangia}

%
%We mention that constant mean curvature hypersurfaces in space forms have been largely studied in literature, we mention the following paper: \cite{dCL} where they prove an Alexandrov-Bernstein Theorem in the hyperbolic space, \cite{KKMS} where they study constant mean curvature surfaces in the hyperbolic space, \cite{LR,rosenberg,NR,GS,GS2,GS3,DSS,ALR,LM,Mon} where they study constant mean curvature hypersurfaces in the hyperbolic space, \cite{GSX} where they study the Plateau problem in the hyperbolic space and the already mentioned paper by Brendle \cite{Brendle}.
%
%
%
%\todo{sfera}

It has been recently shown that Alexandrov's Theorem holds also in a nonlocal setting \cite{Cabre,CFMN}.
We recall that nonlocal minimal surfaces are boundaries of sets $\Omega\subset\mathbb{R}^n$ which are stationary for the $s$-perimeter functional
$$
P_s(\Omega)=\int_\Omega\int_{\Omega^c}\dfrac{dx\, dy}{|x-y|^{n+2s}}\, , \quad \Omega^c=\mathbb{R}^n\setminus\Omega\, ,
$$
with $s\in(0, 1/2)$. 
%If $\Omega$ is an open set with smooth boundary and $A\subset\mathbb{R}^n$ is an open set, %then the condition 
%$$
%\delta P_s(\Omega)=\dfrac{d}{dt}P_s(\Phi_t(\Omega))=0\, , \quad \forall \,  X\in C^\infty_c(A;\mathbb{R}^n)\, ,
%$$
%(where $\Phi_t$ denotes the flux defined by the vector-field $X$) is equivalent to require the vanishing of the nonlocal mean curvature $H^\Omega_s(p)$ of $\Omega$ at every point $p\in\partial A\cap\partial\Omega$. More in general, 
%we say that $H^\Omega_s : A\cap\partial\Omega\rightarrow\mathbb{R}$ is the nonlocal mean curvature of $\partial\Omega$ inside $A$ if
%$$
%\dfrac{d}{dt}P_s(\Phi_t(\Omega))=\int_{\partial\Omega}H^\Omega_s(x)X(x)\cdot\nu_x\, d\mathcal{H}^{n-1}_x \quad \forall \,  X\in C^\infty_c(A;\mathbb{R}^n)\, ,
%$$
%Here $\Phi_t$ denotes the flux defined by the vector-field $X$, $\nu_x$ is the exterior unit normal to $\Omega$ at $x\in\partial\Omega$, and $\mathcal{H}^{n-1}$ denotes the $(n-1)$-Hausdorff measure.
If $\partial\Omega$ is sufficiently smooth, one can 
show that the nonlocal mean curvature of $\partial\Omega$ at a point $p\in\partial\Omega$ is given by
\begin{equation}\label{H_non_loc}
H^\Omega_s(p)=\dfrac{1}{\omega_{n-2}}\int_{\mathbb{R}^n}\dfrac{\chi_{\Omega^c}(x)-\chi_{\Omega}(x)}{|x-p|^{n+2s}}\, dx\, , 
\end{equation}
where $\chi_E$ denotes the characteristic function of a set $E$, $\omega_{n-2}$ is the measure of the $(n - 2)$-dimensional sphere, and the integral is defined in the principal value sense.
By using this notation, the nonlocal version of Alexandrov's Theorem is the following.
\begin{theorem}[\cite{CFMN} and \cite{Cabre}] \label{thm_nonlocal}
	Let $\Omega$ be a bounded open set of class $C^{1,2s}$. The $s$-nonlocal mean curvature $H^{\Omega}_s$ is constant on $\partial \Omega$ if and only if $\partial\Omega$ is a sphere.
\end{theorem}
%\begin{theorem} \label{thm_nonlocal}
%Let $\Omega$ be a bounded domain. If the nonlocal curvature $H_s$ of $\partial \Omega$ is constant, then $\Omega$ is a ball.
%\end{theorem}
Compared to the classical Alexandrov's Theorem \ref{thm_Alex}, we notice that in Theorem \ref{thm_nonlocal} it is not required that $\Omega$ is connected. This is a special feature of the nonlocal setting, where the value of $H_s$ is not a local quantity but it depends on every point of $\Omega$.

We mention that there exist other ways for proving Theorem \ref{thm_Alex}, which are not based on the MMP and involve integral and geometric identities (see \cite{BianchiniCiraoloSalani,Brendle,16BCS,17BCS,MontielRos, Qui Xia,Ros}). Up to now, there is not a corresponding approach in the nonlocal setting and proving Theorem \ref{thm_nonlocal} by using a different approach is an interesting and challenging open problem.

\medskip

Few years after Alexandrov, Serrin used the MMP to study the following problem: let $\Omega$ be a bounded domain in $\mathbb{R}^n$ with boundary of class $C^2$. Suppose there exists a function $u$ satisying the Dirichlet problem 
\begin{equation}\label{eq1}
\Delta u=-1 \textit{ in } \Omega\, , \quad u=0 \textit{ on } \partial\Omega
\end{equation}
together with another boundary condition
\begin{equation}\label{eq2}
\partial_{\nu}u=c \textit{ on } \partial\Omega\, 
\end{equation}
where $\nu$ is the exterior unit normal to $\partial\Omega$ and $c$ is a constant. Must $\Omega$ be a ball? Serrin showed that the answer is affirmative. 

\begin{theorem}[Serrin \cite{Serrin}] \label{thm_Serrin}
	Let $\Omega\subset\mathbb{R}^n$ be a bounded domain whose boundary is of class $C^2$. There exists a solution $u\in C^1(\bar{\Omega})\cap C^2(\Omega)$ to \eqref{eq1} and \eqref{eq2}
	for some constant $c>0$ if and only if $\Omega$ is a ball and $u$ is a radial function.
\end{theorem}

Problems like \eqref{eq1}-\eqref{eq2} are called \emph{overdetermined boundary value problems} since the Dirichlet problem \eqref{eq1} already admits a unique solution; hence condition \eqref{eq2} is an additional requirement and in general the  problem  \eqref{eq1}-\eqref{eq2} may not admit a solution. 

%Thus, the remaining data of the problem, the domain, cannot be given arbitrarly, i.e. there is a requirement also on the domain; that's why the problem is called this phenomenon is called overdetermined and Theorems like Theorem \ref{thm_Serrin} are called \emph{rigidity} phenomena.

We mention that the MMP can be used to prove a more general version of Theorem \ref{thm_Serrin} involving uniformly elliptic quasilinear equation (see \cite{Serrin}). Moreover the MMP has been used to prove an analogue result in space forms (see \cite{Kumaresan-Prajapat, Molzon}) and in a nonlocal setting (see \cite{Jarohs}). 

As for Theorem \ref{thm_Alex}, Theorem \ref{thm_Serrin} can be proven by using different approaches, based on integral identities. This allows to extend Theorem \ref{thm_Serrin} to possibly degenerate quasilinear equations and to fully nonlinear equations (see \cite{AgoMagn,ABennett,BianchiniCiraolo,Henrot5,BiDe,BrandNitsSalTrom,Henrot4,Buttazzo_Kawhol,CianchiSalani,CV_Manuscripta,CHH,Crasta_Fragala,DaLioSir,FarinaKawhol,FrGa,FragalaGazzolaKawhol,GarofaloLewis,HHP,Henrot_Philippin1,Henrot_Philippin2,Henrot_Philippin3,Payne_Sch,Qui Xia bis,Ramm,Ron,Silvestre_Syrakov,Syr,WangXia,Weiberger}) and also for domains with Lipschitz singularities or contained in a convex cone (see \cite{CiraoloRoncoroni,PacellaTralli,Praj}).

For the sake of completeness we mention that these kind of problems arise from practical situations in physics e.g. when one considers the motion of a viscous incompressible fluid moving in straight parallel streamlines through a pipe with planar section $\Omega$ or the torsion of a solid straight bar of given cross section $\Omega$. For these examples, Theorem \ref{thm_Serrin} has the following meaning: ``the tangential stress on the pipe wall is the same at all points of the wall if and only if the pipe has a circular cross section" and ``when a solid straight bar is subject to torsion, the magnitude of the resulting traction which occurs at the surface of the bar is independent of position if and only if the bar
has a circular cross section" (taken from \cite{Serrin}). When the Laplacian is replaced by the mean curvature operator then \eqref{eq1}-\eqref{eq2} describe the shape of a capillary surface in
absence of gravity, adhering to a given plane with constant contact angle. In this example, the analogue of Theorem \ref{thm_Serrin} means that the wetted area on the plane is necessarily spherical (see \cite{Serrin,Wente_appl}). For degenerate elliptic operators, further physical applications may be pointed out. For the $p$-Laplacian \eqref{eq1}-\eqref{eq2} models torsional creep with constant stress on the boundary (see \cite{Kawhol_appl}). For more general operator one has applications in the theory for electrostatic fields (see \cite{Pisani_appl}), we also refer to \cite{Benci, Benci_bis} for more general applications to quantum physics.

\medskip

The last problem that we review in this manuscript is a sort of discrete version of Serrin's overdetermined problem (see \cite{CMS,CMS2,Shagolian}). This problem arises from the study of {\it invariant isothermic surfaces} of a nonlinear non-degenerate {\it fast diffusion} equation (see \cite{MSaihp}), and the main goal is to show that positive solutions of homogeneous Dirichlet boundary value problems or initial-boundary value problems for certain elliptic or, respectively, parabolic equations must be radially symmetric if  {\it just one} of their level surfaces is parallel to $\partial\Omega$ that is,
if the distance of its points from $\partial\Omega$ remains constant.

In order to state the theorem we set up some notations: given a bounded domain $\Omega\subset\mathbb{R}^n$, for $x\in\bar\Omega$ we denote by $d(x)$ the distance of $x$ from $\mathbb{R}^n\setminus\Omega$, that is
$$
d(x)=\min_{y\in\mathbb{R}^n\setminus\Omega}|x-y|\,, \, \, x\in\bar\Omega\, .
$$
For a positive number $\delta$,we define the \emph{parallel surface} to the boundary $\partial\Omega$ of $\Omega$ as 
$$
\Gamma_\delta=\{ x\in\Omega\, :\, d(x)=\delta\}\, .
$$
Moreover, we suppose that there exists a domain $G$ such that 
\begin{equation}\label{palla_interna}
\overline{G}\subset\Omega, \, \partial G\in C^1 \, \text{ satisfying the interior sphere condition, and } \, \partial G=\Gamma_\delta\, ,
\end{equation}
for some $\delta>0$.
\begin{theorem}[\cite{CMS}] \label{thm_Rolly}
	Let $\Omega\subset\mathbb{R}^n$ be a bounded domain and let $G$ satisfy \eqref{palla_interna}. There exists $u\in W^{1,\infty}_0(\Omega)$ solution to \eqref{eq1} such that
	%\begin{equation}\label{eq1_bis}
	%\Delta u=-1 \textit{ in } \Omega, \quad u=0 \textit{ on } \partial\Omega,
	%\end{equation} 
	\begin{equation}\label{eq2_bis}
	u=c \textit{ on } \Gamma_\delta\, ,
	\end{equation}
	for some constant $c>0$, if and only if $\Omega$ is a ball.
\end{theorem}
We mention that Theorem \ref{thm_Rolly} was proved in \cite{CMS} for more general degenerate quasilinear operators as well as for minimizers of not differentiable functionals and for parabolic equations.

We emphasize that in this case the only available proof of this theorem is by using the MMP and, differently from Theorems \ref{thm_Alex} and \ref{thm_Serrin}, it is not available a proof which avoids the MMP and uses instead integral and geometric identities.

\medskip

In the symmetry Theorems \ref{thm_Alex}--\ref{thm_Rolly} the conclusion is the same: the solution to those problems exists if and only if the domain is a ball. The rigidity of these problems is due to the following overdetermined conditions: the mean curvature $H$ is constant in Theorems \ref{thm_Alex} and \ref{thm_nonlocal}, the normal derivative $\partial_{\nu} u$ is constant on $\partial \Omega$ in Theorem \ref{thm_Serrin}, the solution $u$ is constant on $\Gamma_\delta$ in Theorem \ref{thm_Rolly}.

It is natural to investigate the stability of this results: if the overdetermined condition is slightly perturbed, can we say that the domain $\Omega$ is close to a ball? Can we quantify how much the domain is close to a ball?

Stated like that, the assertion is in general false. Indeed, it has been showed in \cite{BNST} for Serrin's Theorem and in \cite{CM} for Alexandrov's Theorem that if the overdetermined condition is close to a constant then the domain $\Omega$ can be close to a bunch of balls connected by small necks (this phenomenon is called \emph{bubbling}). Both \cite{BNST} and \cite{CM} do not invoke the MMP and \emph{perturb} other proofs which use integral identities. 

If some other condition is introduced in order to prevent the bubbling, then the MMP can be studied in a quantitative way to obtain  quantitative information of the proximity of the solution to a single ball. In particular, we will consider the notion of \emph{touching ball condition}.  We observe that the $C^2$-regularity implies a uniform touching ball condition and we denote by $\rho$ the optimal uniform radius, i.e. for any $p\in \partial \Omega$ there exist two balls of radius $\rho$ centered at $c^-\in\Omega$ and $c^+\in\mathbb{R}^{n}\setminus \bar \Omega$ such that $B_\rho(c^-)\subset\Omega$, $B_\rho(c^+)\subset\mathbb{R}^{n+}\setminus \bar \Omega$ and $p\in\partial B_\rho(c^{\pm})$.
%$$
%\begin{tikzpicture}
%\filldraw [black] (2.8,0) circle (1pt) node[right] {\footnotesize $c^+$}
%(2.1,0) circle (1pt);
%\filldraw [black] (1.4,0) circle (1pt) node[left] {\footnotesize $c^-$};
%\node at (-1.7,1.2) {$S$};
%\node at (1.8,0) {\footnotesize $p$};
%\node at (0.5,-0.8) {\footnotesize $B_\rho(c^-)$};
%\node at (2.5,0.95) {\footnotesize $B_\rho(c^+)$};
%\node at (0,0) { $\Omega$};
%\draw(0,0) ellipse (60pt and 40pt);
%\draw(2.8,0) ellipse (20pt and 20pt);
%\draw(1.4,0) ellipse (20pt and 20pt);
%\end{tikzpicture}
%$$

The first quantitative study of the method of moving planes was done by Aftalion Busca and Reichel in \cite{ABR}, where a quantitative version of Serrin's symmetry result was proved. More precisely they proved the following theorem.

\begin{theorem}[\cite{ABR}] \label{thm_ABR}
	Let $\Omega\subset\mathbb{R}^n$ be a bounded domain with $C^{2,\alpha}$ boundary and let $u\in C^2(\overline{\Omega})$ be a positive solution to \eqref{eq1}. There exist two constants $\varepsilon, C>0$ such that the following holds. Assume that   
	\begin{equation}
	\deficit \leq \varepsilon
	\end{equation}
	where 
	$$
	\deficit = ||\partial_\nu u-c||_{C^1(\partial\Omega)} \,,
	$$
	for a constant $c>0$. Then there are two concentric balls $B_r$, $B_R$ such that
	\begin{equation}
	B_r\subseteq\Omega\subseteq B_R\, ,
	\end{equation}
	and 
	\begin{equation}\label{stab_ABR}
	R-r\leq C\left|\log \deficit\right|^{-1/n}\,.
	\end{equation}
	The constants $\varepsilon$ and $C$ depend only on the $C^{2,\alpha}-$regularity of $\Omega$ (and in particular on the radius of the optimal touching ball condition for $\Omega$) and on an upper bound on $\mathrm{diam}(\Omega)$.  
\end{theorem}

We mention that Theorem \ref{thm_ABR} was proved in \cite{ABR} for semilinear equations of the form $\Delta u + f(u)=0$. We notice that the stability estimate \eqref{stab_ABR} is not-optimal, and achieving a sharper stability inequality is an open problem. We also mention that in \cite{BNST} the bubbling phenomenon has been investigated by using an integral identities approach, and in \cite{CirMagnVespri}  the stability estimate \eqref{stab_ABR} has been improved to a power like rate for domains enjoying a further geometric property.

We stress that the assumption that prevent the bubbling in \cite{ABR} is the touching ball condition. This condition will be used also for the other quantitative results reviewed in this manuscript (see Theorems \ref{thm_CMS_stability} and \ref{thm_CV} below) with the exception of the nonlocal Alexandrov's Theorem (see Theorem \ref{thm_CMFG} below).

The following theorem gives a sharp quantitative version of Theorem \ref{thm_Rolly} and it was investigated in \cite{CMS2}. 
\begin{theorem}[\cite{CMS2}] \label{thm_CMS_stability}
	Let $G\subset\mathbb{R}^n$ be a bounded domain with $C^{2,\alpha}$ boundary and set $\Omega= G+B_\delta$ for some $\delta > 0$. Let $u\in C^2(\Omega)\cap C^0(\overline{\Omega})$ be the solution to \eqref{eq1} and let
	$$
	\deficit = \sup_{{\substack{x,y \in \pa G, \\ x\neq y}}}\dfrac{|u(x)-u(y)|}{|x-y|} \,.
	$$
	There exist constants $\varepsilon, C>0$ such that, if
	$$
	\deficit\leq\varepsilon
	$$
	then there are two concentric balls $B_{r}$ and $B_{R}$ such that
	\begin{equation}
	B_r\subseteq\Omega\subseteq B_R\, ,
	\end{equation}
	and 
	\begin{equation}\label{stab_CMS2}
	R-r\leq C\, \deficit\, .
	\end{equation}
	The constants $\varepsilon$ and $C$ only depend on $n$, the $C^{2,\alpha}$-regularity of $\partial G$, the
	diameter of $G$ and $\delta$.
\end{theorem}

The last quantitative results that we review in this manuscript are about Alexandrov soap bubble Theorem (local and nonlocal). We first consider the classical version of Alexandrov's Theorem.  In particular, we consider an $n$-dimensional, $C^2$-regular, connected, closed hypersurface $S$ embedded in $\mathbb{R}^{n+1}$. Given $p\in S$, we denote by $H(p)$ the mean curvature of $S$ at $p$, and we let
$$
\mathrm{osc}(H)=\max_{p\in S}H(p)-\min_{p\in S}H(p)\, .
$$
The quantitative version of Alexandrov's Theorem is the following
\begin{theorem} [\cite{CV_JEMS}]\label{thm_CV}
	Let $S= \partial \Omega$ be an $n$-dimensional, $C^2$-regular, connected, closed hypersurface embedded in $\mathbb{R}^{n+1}$, with $\Omega \subset \mathbb{R}^{n+1}$ a bounded domain satisfying a touching ball condition of radius $\rho$, and let 
	$$
	\deficit = \mathrm{osc}(H) \,.
	$$ 
	There exist constants $\varepsilon, C>0$ such that if
	\begin{equation}
	\deficit \leq\varepsilon\, ,
	\end{equation}
	then there are two concentric balls $B_{r}$ and $B_{R}$ such that
	\begin{equation}
	B_r\subseteq \Omega\subseteq B_R\, ,
	\end{equation}
	and 
	\begin{equation}\label{stab_est}
	R-r\leq C \deficit\, .
	\end{equation}
	The constants $\varepsilon$ and $C$ only depend on $n$ and an upper bound on $\rho^{-1}$ and on $|S|$.
\end{theorem}

Under the assumption that $\Omega$ bounds a convex domain, there exist some results in the
spirit of Theorem \ref{thm_CV} in the literature. In particular, when the domain is an ovaloid, the problem was studied by Koutroufiotis \cite{Kou}, Lang \cite{L} and Moore \cite{Moo}. Other stability results can be found in Schneider \cite{Schneider} and Arnold \cite{Ar}. These results were improved by Kohlmann \cite{Ko} who proved an explicit Hölder type stability in \eqref{stab_est}. In Theorem \ref{thm_CV}, any convexity assumption was done and the rate of stability in \eqref{stab_est} is optimal, as can be proven by a simple calculation for ellipsoids.

We mention that Theorem \ref{thm_CV} can be generalized to hypersurfaces embedded in the Hyperbolic space and in the Hemisphere, see \cite{CV_indiana} and \cite{CRV} respectively. Moreover in \cite{CRV} is shown that the same result holds by replacing the mean curvature with other symmetric elementary functions of the principal curvatures. 
We emphasize that the stability estimate \eqref{stab_est} is optimal. Another optimal stability estimate for proximity to a single sphere was obtained in \cite{KrummelMaggi} by using a different approach. Other quantitative studies regarding the proximity to a single ball can be found in \cite{Fe,magnanini_survey,magnanini_poggesi1,magnanini_poggesi2} where a different deficit is considered.

We notice that Theorem \ref{thm_CV} is not optimal from a qualitative point of view (see \cite{K6} and the reference therein and \cite{Butscher,Butscher_Mazzeo}), since the touching ball condition prevents the possibility of having a bubbling phenomenon. In this direction, Theorem \ref{thm_CV} was improved in \cite{CM} where the bubbling phenomenon was characterized (see also \cite{DMMN} for the anisotropic counterpart).

The quantitative version of Theorem \ref{thm_nonlocal} was investigated in \cite{CFMN} where it is proved that if $H^\Omega_s$ has small Lipschitz constant then $\partial\Omega$ is close to a sphere, with a sharp estimate in terms of the deficit. To state the result we introduce the following uniform distance from being a ball:
$$
\rho(\Omega)=\inf\left\{\dfrac{t-s}{\mathrm{diam}(\Omega)}\, : \, p\in\Omega \, , \,  B_s(p)\subset\Omega\subset B_t(p)  \right\} \, .
$$
\begin{theorem}[\cite{CFMN}] \label{thm_CMFG}
	If $\Omega$ is a bounded open set with $C^{2,\alpha}$ boundary for some $\alpha > 2s$, then there exists a dimensional constant $\hat{C}(n)$ such that
	\begin{equation}\label{rho}
	\rho(\Omega)\leq \hat{C}(n)\dfrac{\mathrm{diam}(\Omega)^{2n+2s+1}}{|\Omega|^2} \deficit\,, 
	\end{equation}
	where
	\begin{equation}\label{deficit_CFMN}
	\deficit=\sup_{{\substack{p,q\in\partial\Omega\, \\ p\neq q}}}\dfrac{|H^\Omega_s(p)-H^\Omega_s(q)|}{|p-q|}\,.
	\end{equation}

\end{theorem}

We emphasize that in Theorem \ref{thm_CMFG} no assumptions on the connectedness of $\Omega$ and on the touching ball condition are done. This is a peculiarity of the nonlocal problem, where every point of $\partial \Omega$ influences the value of the mean curvature at any other point of $\partial \Omega$. Hence, a bubbling result like the one in \cite{CM} is not possible in the nonlocal case. However, it is an open question whether the bubbling can happen for a different type of deficit.

We mention that, under suitable regularity of $\partial \Omega$, from Theorems \ref{thm_CV} and \ref{thm_CMFG} one can prove that $\partial \Omega$ is $C^{1,\alpha}$ close to a single sphere, which means that $\partial \Omega$ can be parametrized as a smooth map of the form $F : \partial B \rightarrow \partial\Omega$ with $\|F - Id\|_{C^{1,\alpha}} \leq C \deficit$. Moreover, in the nonlocal case, we can show also more, that is the $C^{2,\alpha}$ proximity to a single sphere (see \cite[Theorem 1.5]{CFMN}). This result gives an intriguing feature of the nonlocal case, which is the following: if the deficit is small then $\Omega$ is convex (and close to a single sphere).

\medskip

The paper is organized as follows. In Section \ref{The Method of Moving Planes} we introduce some notation which will be used in the rest of the paper. In Section \ref{sect_symmetry} we describe the method of moving planes and recall the proofs of Theorems \ref{thm_Alex}-\ref{thm_Rolly}. In Section \ref{section_almost} we review the quantitative results in Theorems \ref{thm_ABR}--\ref{thm_CMFG}. In Section \ref{sect_open} we state some open problem.

\section{Notation} \label{The Method of Moving Planes}
In this section we introduce some notation which is useful for the application of the MMP. 
Given an arbitrary set $A$ in $\mathbb{R}^n$, a unit vector $\omega\in\mathbb{R}^n$ and a parameter $\lambda\in\mathbb{R}$ we define the following objects:
\begin{equation*}
\begin{array}{lll}
& \pi_\lambda=\{\xi\in\mathbb{R}^n\, :\, \xi\cdot\omega=\lambda\}, \quad & \text{a hyperplane orthogonal to $\omega$} \\
& A^\lambda=\{p\in A\, :\, p\cdot\omega>\lambda\}, & \text{the right-hand cap of $A$}  \\
& \xi^{\lambda} =\xi-2(\xi\cdot\omega-\lambda)\omega, & \text{the reflection of $\xi$ about $\pi_\lambda$}  \\
& A_\lambda =\{p\in\mathbb{R}^n \, : \, p^\lambda\in A^\lambda\}, & \text{the reflected cap about $\pi_\lambda$}, \\
& \hat{A}_\lambda =\{p\in A \, : \, p\cdot\omega<\lambda\}, & \text{the portion of $A$ in the left half-plane}, \\
& \mathcal{M}=\inf\{\lambda\in\mathbb{R}\, : \, A^\lambda=\emptyset\}\, & \text{the extent of $A$ in the direction $\omega$}.
\end{array}	
\end{equation*}
%Lastly, we set $\mathcal{M}$ the extent of $A$ in the direction $\omega$:
%$$
%\mathcal{M}=\inf\{\lambda\in\mathbb{R}\, : \, A^\lambda=\emptyset\}\, ;
%$$ 
%observe that if $\lambda<\mathcal{M}$ is close to $\mathcal{M}$,
%the reflected cap $\Omega_\lambda$ is contained in $\Omega$. Set
%
%.....
The MMP works as follows. Let $\omega$ be a fixed direction and consider the family of hyperplanes $\{\pi_\lambda\}_{\lambda \in \mathbb{R}}$ orthogonal to $\omega$. Since $\Omega$ is bounded, for $\lambda$ very large, $\lambda > \mathcal{M}$, we have that $\pi_\lambda$ does not intersect $\Omega$. We decrease the value of $\lambda$ until $\lambda = \mathcal{M}$, when $\pi_\mathcal{M}$ and $\Omega$ are tangent. Since $\partial \Omega$ is smooth ($C^1$ is enough, see \cite{Fraenkel_libro}), we can decrease $\lambda$ and, at least at the beginning, the reflection $\Omega_\lambda$ of $\Omega^\lambda$ is contained in $\Omega$. We continue in decreasing $\lambda$ until $\Omega^\lambda$ is contained in $\Omega$ and, since $\Omega$ is bounded, we reach a critical value 
$$
\lambda_* = \inf\{ \lambda\in\mathbb{R} \, : \, \Omega_t \subset \Omega \text{ for any } t \in (\lambda,\mathcal{M})\} \,.  
$$
When $\lambda = \lambda_*$, there are two possible critical configurations: 
\begin{itemize}
	\item[$(i)$] $\Omega_{\lambda_*}$ is tangent to $\Omega$ at a point $p$ which is not on $\pi_{\lambda_*}$,
	\item[$(ii)$] $\Omega_{\lambda_*}$ is tangent to $\Omega$ at a point $q$ which is on $\pi_{\lambda_*}$.
\end{itemize}

This is the point where maximum principles enter into play in order to prove symmetry, as we are going to show in the next sections. The rough idea in Theorems \ref{thm_Serrin} and \ref{thm_Rolly} is to compare the solution $u$ to \eqref{eq1} to its reflection $v$, which is defined in $\Omega_{\lambda_*}$ by 
\begin{equation}\label{riflessa}
v(x)=u(x^{\lambda_*})\quad x\in \Omega_{\lambda_*}\,.
\end{equation}
In Theorem \ref{thm_Alex} we will locally compare the surface $S$ to its reflection $S_{\lambda_*}$ by using the fact that $S$ can be locally parametrized over the tangent space to a point by a function which satisfies an elliptic equation (the mean curvature equation). An analogous argument holds for Theorem \ref{thm_nonlocal}, where we can take advantage of the explicit expression of the nonlocal mean curvature \eqref{H_non_loc}.

\begin{rem} \label{rem_MMP_spaceforms}
	As we mentioned in the introduction, the method of the moving planes may be applied in other manifolds, in particular in the hyperbolic space and in the sphere. In these cases, the method needs some more explanation. Here, the notion of hyperplane and direction are replaced by totally geodesic hypersurfaces and geodesic tangent to a fixed direction at a reference point, respectively.
	
	Hence, in the hyperbolic half-space model, the method consists in moving hyperbolic hyperplanes (Euclidean half-spheres or vertical hyperplanes) along a geodesic tangent to a fixed direction at a reference point, say $e_n$. When the critical position is reached, it may be convenient to consider an isometry and to regard the critical hyperplane as a Euclidean vertical hyperplane, so that the reflection is geometrically the Euclidean one.
	
	Analogously, in the hemisphere the method consists in moving half spherical hyperplanes (great circles intersecting the hemisphere) along a geodesic tangent to a fixed direction at a reference point, say the north pole. When the critical position is reached, it may be convenient to consider the stereographic projection and apply a rotation in $\mathbb R^n$ in order to regard the critical hyperplane as a Euclidean vertical hyperplane and, again, the reflection is geometrically the Euclidean one.
	
	%\todo{sup minime}
	
\end{rem}

\section{Symmetry: the qualitative approach} \label{sect_symmetry}
In this section we review the symmetry results of Theorems \ref{thm_Alex}--\ref{thm_Rolly}, which follow by applying the method of moving planes
together with suitable maximum principles.

\medskip 

\subsection{Proof of Theorem \ref{thm_Alex}}\label{Alex_subsect}
We apply the MMP in a fixed direction $\omega$ and, for $\lambda=\lambda_*$, we find the critical positions described in Section \ref{The Method of Moving Planes}. 

If case $(i)$ occurs, then we locally write $S_{\lambda_*}$ and $S$ as graphs of function $u_1$ and $u_0$, respectively, over $B_r\cap T_p S$, where $p$ is the tangency point. Here we denote by $T_p S$ the tangent hyperplane to $S$ at $p$, and we fix a system of coordinates on $T_p S$ such that $p=(0,u_i(0))$ for $i=0,1$. Notice that since $S$ and $S_{\lambda_*}$ are tangent at $p$ then $T_pS$ and $T_p S_{\lambda_*}$ coincide. It is well-known that $u_0$ and $u_1$ satisfy the mean curvature equation (see \cite[Chapter 16]{GT})
$$
{\rm div} \left( \frac{ \nabla u_i(x) }{\sqrt{1+ |\nabla u_i(x)|^2}} \right) = H(x,u(x)) \,, 
$$
$i=0,1$, $x \in B_r\cap T_p S$, for some $r>0$. By construction of the MMP, it is clear that $w=u_1-u_0$ is non-negative, and since the mean curvature $H$ is constant, $w$ satisfies 
$$
Lw=0 \quad \text{in} \, B_r\cap T_p S
$$
for some $r>0$, where $L$ is a linear uniformly elliptic operator (this is possible by choosing $r$ small enough). Since $w(0)=0$, by the strong maximum principle we obtain $w=0$ in $B_r\cap T_p S$, that is, $S$ and $S_{\lambda_*}$ coincide in an open neighborhood of $p$. 

If case $(ii)$ occurs, we locally write $S_{\lambda_*}$ and $S$ as graphs of function $u_1$ and $u_0$ over $T_q S\cap\{x\cdot\omega\leq\lambda_*\}$. As for case $(i)$, we find that there exists $r>0$ such that 
\begin{equation*}
\begin{cases}
L w=0 &\mbox{in } B_r\cap T_q S\cap\{x\cdot\omega<\lambda_*\} , \\ w=0  &\mbox{on } B_r\cap T_q S\cap\{x\cdot\omega=\lambda_*\}.
\end{cases}
\end{equation*}
Since $\nabla w(0)=0$, from Hopf's Lemma we deduce that $w=0$ in $B_r\cap T_q S\cap\{x\cdot\omega\leq\lambda_*\}$. 

Hence, in both cases $(i)$ and $(ii)$ the set of tangency points (that is, those points for which $(i)$ and $(ii)$ occurs) is open. Since it is also closed and non-empty, we must have $S_{\lambda_*}=\hat{S}_{\lambda_*}$, i.e. $S$ is symmetric about the hyperplane $\pi_{\lambda_*}$. Since $\omega$ is arbitrarly, we find that $S$ is symmetric in every direction and we obtain that $S$ is a sphere (see e.g. \cite[Chapter VII, Lemma 2.2]{Hopf}).

%\todo{The geometric version of the Maximum Principle and Hopf Lemma for hypersurfaces.}

%\todo{Differenza tra Serrin, Rolly e Alexandrov... vedi Magnanini Survey}

\medskip

\subsection{Proof of Theorem \ref{thm_Rolly}}
Let $\omega$ be a fixed direction. We apply the method of moving planes to the domain $G$  in the direction $\omega$ and we find the critical positions described for $\lambda=\lambda_*$. A crucial remark is given by the following lemma (see \cite[Lemma 2.8]{CMS}): 
\begin{lemma}
	Let $G$ satisfy assumption \ref{palla_interna}. Then we have
	$$
	\Omega=G+B_\delta(0)=\{ x+y \, : \, x\in G, \, y\in B_\delta(0)\}
	$$
	and, if $G_\lambda\subset G$, then 
	$$
	\Omega_\lambda\subset\Omega \,.
	$$
\end{lemma}
Let $u$ be the solution to \eqref{eq1} and consider the function $v$ defined in $\Omega_{\lambda_*}$ by \eqref{riflessa}. Since $u\geq v$ on $\partial\Omega_{\lambda_*}$, the weak comparison principle yields $u\geq v$ in $\Omega_{\lambda_*}$. Hence the function $w = u-v$ satisfies 
$$
\Delta w = 0 \quad \text{ and } \quad w = u-v \geq 0 \quad \text{ in }  \Omega_{\lambda_*} \,.
$$
If case $(i)$ of the MMP occurs, we apply the strong maximum principle to $w$ in $\Omega_{\lambda_*}$ and obtain that $w>0$ in $\Omega_{\lambda_*}$, since $w \not \equiv 0 $ on $\Gamma_\delta\cap\Omega_{\lambda_*}$. This is a contradiction, since $p$ belongs both to $\Omega_{\lambda_*}$ and $\Gamma_\delta\cap\Gamma_\delta^{\lambda_*}$, and hence $u(p)=v(p)$, i.e. $w(p)=0$.

Now, let us consider case $(ii)$ of the MMP. Notice that $\omega$ belongs to the tangent hyperplane to $\Gamma_\delta$ at $q$. By applying Hopf's lemma in $\Omega_{\lambda_*}$ we get 
\begin{equation*}
\partial_\omega w(q)< 0 \, .
\end{equation*}
On the other hand, since $\Gamma_\delta$ is a level surface of $u$ and $u$ is differentiable at $q$, we must have that
\begin{equation*}
\partial_\omega w(q)=0\, ;
\end{equation*} 
this gives the desired contradiction.

Hence, we have proved that for any direction we have that $\Omega = \Omega_{\lambda_*}$, i.e. $\Omega$ is symmetric in any direction $\Omega$. This implies that $\Omega$ is a ball and $u$ is radially symmetric.

\medskip

\subsection{Proof of Theorem \ref{thm_Serrin}}
As usual, we apply the MMP in a fixed direction $\omega$ and we stop at a critical position for $\lambda=\lambda_*$. 

In order to prove that $\Omega$ is symmetric with respect to the hyperplane $\pi_{\lambda_*}$ we consider the function $w = u-v$ in $\Omega_{\lambda_*}$, where $v$ defined  by \eqref{riflessa}. It is clear that $w$ satisfies
\begin{equation*}
\begin{cases}
\Delta w=0 &\mbox{in } \Omega_{\lambda_*} , \\ w=0  &\mbox{on } \partial \Omega_{\lambda_*}\cap\pi_{\lambda_*}, \\ w\geq 0 &\mbox{on } \partial \Omega_{\lambda_*}\setminus\pi_{\lambda_*}.
\end{cases}
\end{equation*}
At this point the strong maximum principle gives either 
\begin{equation}\label{assurdo}
w>0 \quad \text{in} \, \Omega_{\lambda_*}
\end{equation}
or $w \equiv 0$ in $\Omega_{\lambda_*}$. The latter case would imply that $\Omega$ is symmetric about $\pi_{\lambda_*}$. Hence, we assume that $w>0$ in $\Omega_{\lambda_*}$.

Assume that case $(i)$ occurs, that is $\Omega_{\lambda_*}$ is internally tangent to $\partial\Omega$ at a point $p$ not belonging to $\pi_{\lambda_*}$ and assume by contradiction that \eqref{assurdo} holds true. Then Hopf's Lemma ensures that 
\begin{equation*}
\partial_\nu w(p)<0\, ,
\end{equation*} 
but this contradicts the fact that  $\partial_\nu w(p) = 0$, since
\begin{equation}
\partial_\nu u(p)=\partial_\nu v(p)=c \, .
\end{equation} 
We conclude that \eqref{assurdo} cannot occur in case $(i)$.

Case $(ii)$ is much more complicated since Hopf's Lemma cannot apply. The proof makes use of a refinement of the maximum principle, see Lemma \ref{Lemma_Serrin} below (for its proof see \cite{Serrin}). The goal is to prove that $w$ has a second order zero in $q$. To do this we fix a coordinate system with the
origin at $q$, the $x_n$ axis in the direction of the inward normal to $\partial\Omega$ at $q$, and the $x_1$ axis in the direction of $\nu$, that is normal to $\pi_m$. In this coordinate system the boundary of $\Omega$ is locally given by
$$
x_n=\phi(x_1,\dots,x_{n-1})\, \quad \phi\in C^2\,.
$$
Since $u\in C^2$, by the boundary conditions $u=0$ and $\partial_\nu u=c$ on $\partial\Omega$ and differentiating twice, one obtains that 
\begin{equation}\label{u=0_tris}
\partial^2_{x_i x_j}u+c\partial^2_{x_i x_j}\phi=0\,, \quad i, j = 1, \dots , n - 1\,,
\end{equation}
\begin{equation}\label{italy_loves_food}
\partial^2_{x_n x_i}u(q)=0\,,\quad i = 1, \dots , n - 1 \,,
\end{equation}
and
\begin{equation}\label{Laplaciano}
\partial^2_{x_n x_n}u(q)=-\sum_{i=1}^{n-1}\partial^2_{x_i}u(q)-1=c\Delta\phi(q)-1 \,.
\end{equation}
By construction $\Omega_{\lambda_*}\subseteq\Omega$ and $\partial^2_{x_1 x_j}u(q)=0$ for $j = 2, \dots , n-1$, because $\partial_{x_1}\phi$ has an extremum point at $q$ with respect to all but the first coordinates directions. 

Since 
\begin{equation}
v(x_1,x_2,\dots,x_n)=u(-x_1,x_2,\dots,x_n)
\end{equation}
by \eqref{u=0_tris}, \eqref{italy_loves_food} and the last remark we have that all the first and second derivatives of $u$ and $v$ coincide at $q$, and hence 
\begin{equation} \label{wzero}
\nabla w (q)=0 \quad \text{ and } \quad \nabla^2 w (q) =0 \,.
\end{equation} 
On the other hand,  $w$ satisfies
\begin{equation}
\Delta w=0 \quad \text{and} \quad w>0 \quad  \text{in} \quad \Omega_{\lambda_*}
\end{equation}
and $w(q)=0$. The contradiction is obtained thanks to the so called Serrin's corner lemma (see \cite{Serrin} for the proof).

\begin{lemma}\label{Lemma_Serrin}
	Let $\Omega\subset\mathbb{R}^n$ be a $C^2$ bounded domain of $\mathbb{R}^n$ and let $\xi$ be a direction such that $\xi\cdot\nu= 0$ in $y\in\partial\Omega$. Let $H(\nu)$ be an open halfspace with unit outer normal $\nu$, $\Omega(\nu)=\Omega\cap H(\nu)$ and let $w\in C^2(\overline{\Omega(\nu)})$ satisfy
	$$
	\Delta w\leq 0 \quad  w\geq 0 \quad \text{in} \quad \Omega(\nu)\, , \quad \text{and} \quad w(y)=0\, .
	$$ 
	If $\theta$ is a direction in $y$ entering $\Omega(\nu)$ such that $\theta\cdot\nu\neq 0$, then either
	$$
	\partial_\theta w(y)>0 \quad \text{ or } \quad \partial^2_{\theta\theta} w(y)>0\, ,
	$$
	unless $w\equiv 0$.
\end{lemma}

Indeed, let $\theta$ be any direction not parallel to $\nu$. Lemma \ref{Lemma_Serrin}  ensures that either
$$
\partial_\theta w(q)>0 \quad \text{ or } \quad \partial^2_{\theta\theta} w(q)>0\, ,
$$
which is a contradiction on account of \eqref{wzero}.

Hence, we have proved that $\Omega$ is symmetric with respect to the hyperplane $\pi_{\lambda_*}$ and Theorem \ref{thm_Serrin} follows since $\Omega$ is symmetric with respect to any direction $\omega$. Moreover, by construction, $\Omega$ is also simply connected, and then it is a ball and $u$ must be radial.

\subsection{Proof of Theorem \ref{thm_nonlocal}} 
In this subsection we give a sketch of the proof of the nonlocal version of Alexandrov's Theorem, for the details we refer to  \cite{CFMN} and \cite{Cabre}.

We apply the MMP in a direction $\omega$ until we reach the critical position $\lambda=\lambda_*$, and we prove the following inequality: 
\begin{equation}\label{ineq_CFMN}
\int_{\Omega\bigtriangleup\Omega_{\lambda_*}} d(x,\pi_{\lambda_*})\, dx\leq 0\, ;
\end{equation}
where $E\bigtriangleup F$ denotes the symmetric difference of the two sets, that is $E\bigtriangleup F=(E\setminus F)\cup(F\setminus E)$.

In order to prove inequality \eqref{ineq_CFMN}, we can assume without loss of generality that $\omega = e_1$. 

Assume that case $(i)$ of the MMP occurs. Hence there exists $p\in\partial\hat{\Omega}_{\lambda_*}\setminus\pi_{\lambda_*}$ such that $p\in\partial\Omega\cap\partial\Omega_{\lambda_*}$. Then
\begin{equation}\label{3.7}
H^\Omega_s(p)-H^\Omega_s(p')=\dfrac{2}{\omega_{n-2}}\int_{\Omega_{\lambda_*}\setminus\Omega}\left(\dfrac{1}{|x-p|^{n+2s}}-\dfrac{1}{|x^{\lambda_*}-p|^{n+2s}}\right)\, dx\, ,
\end{equation}
where all the integrals are intended in the principal value sense. Since $x^{\lambda_*}=(2\lambda_*-x_1,x_2,\dots,x_n)$ and by the convexity of the function $f(t)=(1+t)^{(n+2s)/2}-1$, we get that if $x\in\Omega_{\lambda_*}$ then
\begin{equation}\label{3.9}
\begin{aligned}
\dfrac{1}{|x-p|^{n+2s}}-\dfrac{1}{|x^{\lambda_*}-p|^{n+2s}}\geq&\dfrac{2(n+2s)(x_1-\lambda_*)(p_1-\lambda_*)}{|x^{\lambda_*}-p|^{n+2s}|x-p|^2} \\ \geq&\dfrac{2(n+2s)(x_1-\lambda_*)(p_1-\lambda_*)}{\mathrm{diam}(\Omega)^{n+2s+2}}\, ,
\end{aligned}
\end{equation}
where we also used the fact that, by construction, $p^{\lambda_*}\in\partial\Omega$ and therefore 
$$
|x-p|=|x^{\lambda_*}-p^{\lambda_*}|\leq\mathrm{diam}(\Omega) \,.
$$
From \eqref{3.7} and \eqref{3.9} we get \eqref{ineq_CFMN}.

Now, assume that case $(ii)$ of the MMP occurs, i.e. $\pi_{\lambda_*}$ is orthogonal to $\partial\Omega$ at some point $q\in\partial\Omega\cap\pi_{\lambda_*}$. In this case we consider the following quantity: $\nabla H^{\Omega}_s(q)\cdotp e_1$.  Via an approximation argument one can show that 
\begin{equation*}
\nabla H^{\Omega}_s(q)\cdotp e_1=-\dfrac{2(n+2s)}{\omega_{n-2}}\int_{\Omega\setminus\Omega_{\lambda_*}}\dfrac{(x-q)\cdotp e_1}{|x-q|^{n+2s+2}}\, dx\, .
\end{equation*}
Since $\nabla H^{\Omega}_s(q)\cdotp e_1=0$ and $q\in\pi_{\lambda_*}$, we conclude that
\begin{equation*}
\begin{aligned}
0&=\dfrac{2(n+2s)}{\omega_{n-2}}\int_{\Omega\setminus\Omega_{\lambda_*}}\dfrac{(x-q)\cdotp e_1}{|x-q|^{n+2s+2}}\, dx \\
&\geq  \dfrac{2(n+2s)}{\mathrm{diam}(\Omega)^{n+2s+2}\omega_{n-2}}\int_{\Omega\setminus\Omega_{\lambda_*}}|x_1-\lambda_*|\, dx \\
&=\dfrac{2(n+2s)}{\mathrm{diam}(\Omega)^{n+2s+2}\omega_{n-2}}\int_{\Omega_{\lambda_*}\bigtriangleup\Omega}|x_1-\lambda_*|\, dx\, ;
\end{aligned}
\end{equation*}
which completes the proof of \eqref{ineq_CFMN}.

It is clear that \eqref{ineq_CFMN} implies that $|\Omega\bigtriangleup\Omega_{\lambda_*}|=0$. Since the direction is arbitrary,  the argument works in every direction and we conclude that $\Omega$ is symmetric in any direction and than that $\Omega$ is a ball.

\section{Almost symmetry: the quantitative approach} \label{section_almost}
In this section we show how the MMP can be studied from a quantitative point of view in order to obtain the quantitative stability results in Theorems \ref{thm_ABR}--\ref{thm_CMFG} in terms of the deficit $\deficit$ (notice that the deficit changes according to the type of problem that we consider). There is a partial common strategy in the proofs of these theorems, which we explain below. Further details will be given in the subsections below, according to the type of problem which we consider. 

\medskip 

The starting point for studying the MMP from a quantitative point of view is to find an approximate center of symmetry $O$ for the problem. This is achieved by applying the method of moving planes in $n$ orthogonal directions, say $e_1,\ldots,e_n$, and considering the intersection between the corresponding critical planes. More precisely, by applying the MMP in the direction $e_i$, $i=1,\ldots,n$, we obtain the critical hyperplanes $\pi_{\lambda_i}$, $i=1,\ldots,n$, and we define
$$
O = \bigcap_{i=1}^n \pi_{\lambda_i} \,.
$$ 
Notice that, up to a translation we may assume that $O$ is the origin of $\mathbb{R}^n$.

In order to prove that $O$ is an approximate center of symmetry, we need to prove a theorem of approximate symmetry in one direction and to provide a quantitative estimate:

\medskip

\begin{itemize}
	
	\item[\emph{Step 1:}] Quantitative approximate symmetry in any direction in terms of $\deficit$.  
	
\end{itemize}

\medskip

This is the main point of the quantitative proofs and it is used not only for proving that $O$ is an approximate center of symmetry, but it is also invoked several times at other points of the proof of the main results. This step differs in the proofs of Theorems \ref{thm_ABR}--\ref{thm_CMFG}, and its description will be given later in more details in subsections \ref{subsec_quant_1}--\ref{subsec_quant_4} below, according to the type of problem that we consider.

Without entering now in the details of how to prove Step 1, we give an idea of the framework needed to achieve the quantitative estimates. We mention that the framework described here was used in \cite{CV_JEMS} and \cite{CFMN} and differs in some step from the ones adopted in \cite{ABR} and \cite{CMS}.

Let $\omega$ be a fixed direction and let $\lambda_*$ be as in Section \ref{The Method of Moving Planes}. The approximate symmetry in one direction quantifies the symmetric difference between $\Omega$ and $\Omega_{\lambda_*} \cup \Omega^{\lambda_*}$ (i.e. the maximal cap and its reflection about the critical hyperplane $\pi_{\lambda_*}$). More precisely, we denote by $\Sigma$ the connected component of $\Omega_{\lambda_*}$ containing the touching point between $\Omega$ and $\Omega_{\lambda_*}$ (of course, a priori this may be not unique since it is possible that there are more touching points between the two sets). We also set $\Sigma'$ as the reflection of $\Sigma$ about $\pi_{\lambda_*}$ (notice that $\Sigma'$ is a connected component of $\Omega_{\lambda_*}$). The approximate symmetry in one direction implies an estimate of the following form:
\begin{equation} \label{stima1}
|\Omega \triangle (\Sigma \cup  \Sigma') | \leq C \deficit \,.
\end{equation}
Since the reflection $\mathcal{R}^O$ with respect to $O$ can be seen as the composition of $n$-orthogonal reflections, from  \eqref{stima1} we find that it satisfies
\begin{equation} \label{stima2}
|\Omega \triangle (\Sigma \cup  \mathcal{R}^O(\Sigma)) | \leq C \deficit \,.
\end{equation}

Once we have \eqref{stima1} we can prove the following results.

\medskip

\begin{itemize}
	\item[\emph{Step 2:}] For any direction $\omega$, the corresponding critical hyperplane $\pi_{\lambda_*}$ stops close to $O$. In particular
	\begin{equation} \label{distOpi}
	{\rm dist}(O, \pi_{\lambda_*}) \leq C \deficit \,,
	\end{equation}
	where the constant $C$ does not depend on the direction. 
	
	This step is achieved by combining the almost symmetry of $\Omega$ with respect to $O$ and to $\pi_{\lambda_*}$. Indeed, by assuming without loss of generality that $\pi_{\lambda_*}=\{x \cdot \omega = \lambda_* \}$, with $\lambda_* > 0$, it is possible to show that $\Omega$ has small volume in the strip $\{|x\cdot \omega| \leq \lambda_*\}$.  By iterating this argument (reflecting with respect to $O$ and to $\pi_{\lambda_*}$) one can find the bound \eqref{distOpi} (see for instance \cite[Lemma 4.1]{CFMN}).

	\medskip
	
	\item[\emph{Step 3:}] Once \eqref{distOpi} is proved, we define 
	$$
	r= \sup \{s>0 :\ B_s(O) \subset \Omega \} \quad \text{ and } \quad R=\inf \{s>0:\  \Omega \subset B_s(O) \}
	$$
	and prove the bound
	\begin{equation}\label{Rmenor}
	R - r \leq C \deficit \,.
	\end{equation}
	Moreover, $\partial \Omega$ is a Lipschitz graph over the sphere $\partial B_r(O)$. 
	
	Estimate \eqref{Rmenor} is obtained by applying the MMP again and using \eqref{distOpi}. More precisely, by denoting by $x$ and $y$ two points on $\partial \Omega \cap B_r(0)$ and $\partial \Omega \cap B_R(0)$, respectively, we apply the MMP in the direction
	$$
	\omega = \frac{y-x}{|y-x|} \,,
	$$
	and show that $R-r \leq 2 {\rm dist}{(O, \pi_\omega)}$, and from \eqref{distOpi} we obtain \eqref{Rmenor}.
	
	The estimate \eqref{Rmenor} completes the first part of the assertion of the quantitative Theorems \ref{thm_ABR}--\ref{thm_CMFG} which says that $\Omega$ is close to a ball in $L^\infty$ norm. 
	
	As a by product of this proof, one can prove that $\partial \Omega$ is a graph over $\partial B_r(O)$, i.e. there exists a map $\Phi: \partial B_r(O) \to \mathbb{R}^n$ of class $C^{0,1}$ and such that $\Phi(\partial B_r(O))= \partial \Omega$.
	
\end{itemize}

\medskip

\noindent In the case of Alexandrov's Theorem (classical and nonlocal), we can do a further step:

\medskip

\begin{itemize}
	
	\item[\emph{Step 5:}] $\partial \Omega$ is a $C^{1,\alpha}$ small perturbation of the identity map over the sphere. More precisely, we can show that there exists a map $\Phi: \partial B_r(O) \to \mathbb{R}^n$ such that $\Phi(\partial B_r(O))= \partial \Omega$ and 
	\begin{equation} \label{stimaC1}
	\|\Phi - Id\|_{C^{1,\alpha}(\partial B_r(O))} \leq C \deficit \,.
	\end{equation}
	This step is achieved by using regularity theory: from Step 4 the map $\Phi$ is Lipschitz continuous and we can apply elliptic regularity theory to obtain \eqref{stimaC1}.  
	
\end{itemize}

\medskip

In the following subsections, we give some detail regarding how to prove Step 1. i.e. the approximate symmetry in one direction.

\subsection{Serrin's overdetermined type problems} \label{subsec_quant_1}
The first quantitative study of the method of moving planes was done in \cite{ABR} for Serrin's overdetermined problem. The main idea is to replace qualitative tools like the maximum principle with quantitative tools like Harnack's inequality, and to study Hopf Lemma and Serrin's corner lemma from a quantitative point of view. The quantitative bounds are given in terms of the deficit
$$
\deficit := \| \partial_{\nu} u - c\|_{C^1(\partial \Omega)} \,.
$$
The first goal is to give quantitative bounds on the difference function $w=u-v$ between the solution $u$ and the reflection $v$ (see subsection \ref{The Method of Moving Planes} for the definition of $u$ and $v$). Such bounds are given in a subset $D_\delta$ of the connected component $D$ of $\Omega_{\lambda^\ast}$ which contains the tangency point $p$. The set $D_\delta$ is defined as the set of points of $D$ which are far from $\partial D$ more than $\delta$.

More precisely, the Hopf Lemma used in case $(i)$ is replaced by the following inequality (see \cite[Proposition 2]{ABR})
\begin{equation} \label{Hopf_quant}
w(y) \leq C^{1/\delta} |\partial_{\nu} w (p)|
\end{equation}
for any $y \in D_\delta$, where $p$ is the tangency point in case $(i)$ of the MMP, and $C$ depends only on the dimension $n$, $|\Omega|$ and on $\rho$. 

Analogously, the quantitative version of Serrin's corner lemma is given by
\begin{equation} \label{Serrin_corner_quant}
w(y) \leq \frac{C^{1/\delta}}{(s \cdot \omega) |s \cdot \nu(q)|}  |\partial^2_{ss} w (q) |
\end{equation}
for any $y \in D_\delta$, where $q$ is the tangency point in case $(ii)$ of the MMP and $s$ in any direction that enters $D$ nontangentially at $q$. 

Estimates \eqref{Hopf_quant} and \eqref{Serrin_corner_quant} are obtained by using the usual barrier functions used in Hopf's boundary point Lemma and Serrin's corner Lemma, respectively, and then by propagating the obtained information by using Harnack's inequality (this is one of the reasons why we need to restrict the domain $D$ to a smaller domain $D_\delta$).

By using \eqref{Hopf_quant} and \eqref{Serrin_corner_quant}, one can prove that the parallel subset $\Sigma_\delta$ of the maximal cap $\Sigma$ in the direction $\omega$ has a connected component $\tilde \Sigma_\delta$ such that 
\begin{equation} \label{prop1ABR}
\|w\|_{L^\infty(\tilde \Sigma_\delta)} \leq  C^{({\rm diam}(\Omega)/\delta)^n} \deficit  \,,
\end{equation}
where $C$ does not depend on $\delta$.
This estimate says that the solution $u$ and its reflection $v$ are close in $L^\infty$ norm. To obtain information on the almost symmetry of $\Omega$ in the direction $\omega$, we need to do some more step. In particular, we can use that $u$ (and hence $v$) grows linearly from the boundary, which implies that
$$
{\rm dist}(x, \partial \Omega) \leq C u(x)  = C(v(x) + w(x))
$$ 
for any $x \in \tilde \Sigma_\delta$; since $v(x)$ grows linearly from the boundary and from \eqref{prop1ABR} we conclude that 
\begin{equation} \label{stimaABR}
{\rm dist}(x, \partial \Omega) \leq C u(x) \leq C \delta + C_*^{({\rm diam}(\Omega)/\delta)^n} \deficit 
\end{equation}
for any $x \in \partial \Sigma$. 

The final goal is to show that there exists a symmetric set (with respect to $\pi_\omega$) which approximates well $\Omega$. We define this set $X_\sigma$ as 
$$
X_\sigma = (\Sigma \cup \Sigma')_\sigma := \{x \in \Sigma \cup \Sigma':\ {\rm dist}(x \partial(\Sigma \cup \Sigma')>\sigma) \} \,.
$$

From \eqref{stimaABR} one can conclude that, if $\sigma$ and $\delta$ are small compared to $\rho$ and if
\begin{equation} \label{sigma2}
\sigma >  C\left( \delta + C_*^{({\rm diam}(\Omega)/\delta)^n} \deficit \right) \,,
\end{equation}
then
\begin{equation} \label{XsigmaOmega}
\Omega_\delta \subset X_\sigma \subset \Omega \,.
\end{equation}
By choosing $\sigma $ and $\delta$ appropriately, from \eqref{sigma2} we obtain that \eqref{XsigmaOmega} holds provided that $\delta$ is small compared to $\rho$ and 
\begin{equation} \label{abrsigma}
\delta < \sigma < C | \log ( \deficit ) |^{-1/n} \,.
\end{equation}

\noindent As a consequence of \eqref{XsigmaOmega} and \eqref{abrsigma} one obtains the following statement of the approximate symmetry in one direction:

\begin{prop} \label{approx1dABR}
	For any point $x \in \partial \Omega$ there exists a point $y \in \partial \Omega$ such that
	$$
	|x^{\lambda_*}-y| < C | \log ( \deficit ) |^{-1/n} \,,
	$$
	where $x^{\lambda_*}$ is the reflection of $x$ about $\pi_{\lambda_*}$ and $C$ depends only on $\rho,n$ and $|\Omega|$.
\end{prop}

\subsection{Discrete type Serrin's problem}
The approximate symmetry in one direction needed in Theorem \ref{thm_Rolly} is similar to the one for Serrin's overdetermined problem. With respect to Serrin's problem, here we can work inside the domain and we are able to obtain sharp stability estimates in terms of the deficit
\begin{equation} \label{deficit_rolly}
\deficit := \sup_{\substack{x,y \in \partial G, \\ x \neq y}} \frac{|u(x)-u(y)|}{|x-y|} \,.
\end{equation}
In this subsection we give the main ideas for proving the approximate symmetry in one direction for Theorem \ref{thm_CMS_stability} and we refer to \cite{CMS2} for the details. 

The big advantage in this problem is that we can apply the MMP procedure to $G$, which is in the interior of $\Omega$. Moreover, Serrin's corner lemma is not needed, and we just have to replace the strong maximum principle by the Harnack's inequality (in case $(i)$) and to write a quantitative form of Hopf's boundary point lemma (for case $(ii)$). 

More precisely, assume that case $(i)$ in the MMP occurs and that the tangency point $p \in \partial G$ is far from $\pi_{\lambda_*}$ more than $\rho/4$. The point $p$ determines a connected component $\Sigma_G$ of the maximal cap $G^{\lambda_*}$ in the direction $\omega$. Since $p$ is an interior point of $\Omega^{\lambda_*}$, we can apply Harnack's inequality to $w$ and write
\begin{equation} \label{harnack1rolly}
\sup_{B_r(p)} w \leq C \inf_{B_r(p)} w \leq C w(p) \,,
\end{equation}
provided that $B_{2r}(p) \subset \Omega$. Since $p$ is a tangency point, $p,p' \in \partial G$ with $|p-p'| \geq \rho/2$ and 
$$
w(p)=u(p)-v(p') \leq \frac{2}{\rho} \deficit 
$$
and from \eqref{harnack1rolly} we obtain
\begin{equation} \label{harnack1rolly2}
\sup_{B_r(p)} w \leq C\deficit \,.
\end{equation}
By using a suitable chain of balls, we can iterate the use of Harnack's inequality and obtain
\begin{equation} \label{harnack1rolly3}
\sup \{w(x):\ x \in \Sigma  \textmd{ and } {\rm dist}(x,\pi_{\lambda_*}>r \} \leq C\deficit \,,
\end{equation}
where $\Sigma$ is the connected component of $G$ such that $p \in \partial \Sigma$. By using Carleson estimates \cite[Theorem 1.3]{BCN} we can extend \eqref{harnack1rolly3} to the whole $\Sigma$ and obtain
\begin{equation} \label{harnack_final_rolly} 
\sup_{\Sigma} w \leq C\deficit \,.
\end{equation}

If case $(ii)$ of the MMP occurs, instead of  \eqref{harnack1rolly} we have to do a quantitative study of Hopf boundary point lemma. More precisely, if $q \in \partial G \cap \pi_{\lambda_*}$ is a tangency point, then one can prove that 
$$
\sup_{B_{r/2}(p-r\omega)} w \leq C \partial_{\omega} w(q) \,.
$$ 
Since 
$$
\partial_{\omega}u(q) = -\partial_{\omega} v(q) \,.
$$
from the definition of the deficit \eqref{deficit_rolly}, it follows that 
\begin{equation*}  
\sup_{B_{r/2}(q-r\omega)} w \leq \deficit \,
\end{equation*}
and \eqref{harnack_final_rolly} follows by iterating Harnack's inequality and arguing as before.

The final step towards the approximate symmetry in one direction is to prove that if $\deficit$ is small then 
$$
G_s \subset \Sigma \cup \Sigma' \subset G 
$$
for any $s \in (C \deficit , \rho/2)$. This implies that $ \Sigma \cup \Sigma'$ approximates well $G$ (see \cite[Theorem 3.5]{CMS2}). Finally, the approximate symmetry in one direction follows easily.

\begin{prop} \label{approx1d_rolly}
	For any point $x \in \partial G$ there exists a point $y \in \partial G$ such that
	$$
	|x^{\lambda^\ast}-y| < C  \deficit   \,,
	$$
	where $x^{\lambda^\ast}$ is the reflection of $x$ about $\pi_{\lambda^\ast}$ and $C$ depends only on $\rho,n$ and $|G|$.
\end{prop}

\subsection{Alexandrov soap bubble Theorem}
In this subsection we show how to obtain the approximate symmetry in one direction for Alexandrov's Theorem. We use the following deficit
$$
\deficit := {\rm osc} (H) \,.
$$
We first consider the case in which $p \in \partial \Omega$ is a tangency point which is far from $\pi_{\lambda_*}$ more than some fixed quantity $\delta$ which is a small ratio of $\rho$. As done for the symmetry case, we \emph{locally} write $\partial \Omega$ and its reflection about $\pi_{\lambda_*}$ as graph of function on the tangent space to $p$. A crucial observation is that, since $\Omega$ satisfies a touching ball condition of radius $\rho$, then we can quantify the word  \emph{locally} and say that $\partial \Omega$ is locally the graph of function on a ball of radius $\rho$ in the tangent space.  

By arguing as in subsection \ref{Alex_subsect}, we have two functions $u,\hat u : B_\rho(O) \cap T_p(\partial \Omega) \to \mathbb{R}$ which satisfy 
$$
L(u-v)(x) = H(x,u(x)) - H(x,v(x)) \quad \text{and} \quad u - \hat u  \geq 0 \,.
$$
Since the right hand side can be bounded in terms of $\deficit$, Harnack's inequality (see \cite[Theorems 8.17 and 8.18]{GT}) gives that
$$
\sup_{B_\delta(0)} (u-v) \leq C \left( \inf_{B_\delta(0)} (u-v) + \deficit \right)
$$
and since $u(0)=v(0)$ we obtain 
$$
\sup_{B_\delta(0)} (u-v) \leq C \deficit \,.
$$
By using interior regularity estimates (see \cite[Theorem 8.32]{GT}) we obtain that 
\begin{equation} \label{harn1}
\|u-v\|_{C^1(B_\delta(0))} \leq C \deficit \, .
\end{equation}
Estimate $\eqref{harn1}$ says that in a neighborhood of fixed size of $p$, the two surfaces $\Sigma$ and $\hat \Sigma$ are close in $C^1$ norm, so that we can conclude that:

\emph{for any $q \in \mathcal{U}_r(q)$, where $\mathcal{U}_r(q)$ is a neighbourhood of $q$ in $\Sigma$ of fixed size $r$ (see \cite[Formula 2.1]{CV_JEMS} for the precise definition), there exists $\hat q \in \hat\Sigma$ such that }
\begin{equation} \label{step1alex}
|q -\hat q| + |\nu_q - \nu_{\hat q}| \leq C \deficit \,.
\end{equation}
In order to propagate the smallness estimate \eqref{step1alex}  to the whole $\Sigma$, we need to construct a chain of balls contained in $\Sigma$ and iterate the use of Harnack inequality. This iteration is not trivial since at any step we have to use a different system of coordinates and new functions $u$ and $\hat u$. However, the fact that the normals are $C^0$-close guarantees that the change of coordinate system introduce a small error of size $\deficit$. Moreover, the extension of \eqref{step1alex} to the points close to $\pi_{\lambda_*}$ is difficult and one has to argue carefully. We do not enter in the details and refer to \cite[Case 1 pages 283-287]{CV_JEMS}. Moreover, the proof of the approximate symmetry in one direction consists of other three cases which have to be considered separately, according to how much the touching point is close to $\pi_{\lambda_*}$ (see \cite[Case 2,3,4 pages 288-292]{CV_JEMS}).
The final estimate that we prove is the approximate symmetry in one direction. 
\begin{prop}  \label{approx1dCV}
	For any $q \in \Sigma$ there exists $\hat q \in \hat\Sigma$ such that 
	\begin{equation} \label{step1alexbis}
	|q -\hat q| + |\nu_q - \nu_{\hat q}| \leq C \deficit \,.
	\end{equation}
\end{prop}

\subsection{Nonlocal Alexandrov Theorem} \label{subsec_quant_4}
In the nonlocal case one can obtain directly the estimate \eqref{stima1} which, in the other cases, we obtain as a corollary of Propositions \ref{approx1dABR}, \ref{approx1d_rolly} and \ref{approx1dCV}. Indeed, one can explore the proof of Theorem \ref{thm_nonlocal} from a quantitative point of view and obtain the bound
\begin{equation*}
|\Omega \triangle \Omega'| \leq C {\rm diam}(\Omega)^{n+s+(1/2)} \sqrt{\deficit} \,,
\end{equation*} 
which gives the approximate symmetry in one direction. We notice that in this case the approximate symmetry in one direction does not give the sharp exponent (which is linear). However, the sharp stability estimate in Therorem \ref{thm_CMFG} is obtained later by using again a quantitative analysis of the method of moving planes (see \cite[Step 1 in the proof of Theorem 1.2]{CFMN}). Indeed, we set
\begin{equation} \label{def:ri_re}
r=\min_{x\in \partial \Omega} |x|\,,\qquad R=\max_{x\in \partial\Omega} |x|\,,
\end{equation}
and we give an upper bound on $R-r$.

Let $x,y \in \partial\Omega$ be such that $|x|=r$ and $|y|=R$. If $x=y$ then the conclusion is trivial. Hence, we assume that $x \neq y$, consider the direction
\[
e=\frac{y-x}{|y-x|}\,,
\]
and denote by $\pi_{\ell_e}$  the corresponding critical hyperplane. One can prove that $y$ is closer than $x$ to the critical hyperplane $\pi_{\ell_e}$, i.e.,
\begin{equation} \label{dist_x_dist_y}
{\rm dist}(x,\pi_{\ell_e}) \geq {\rm dist}(y,\pi_{\ell_e})\,.
\end{equation}
From \eqref{dist_x_dist_y} and since $e$ is parallel to $y-x$ we get
\begin{equation}
\label{eq:Rr lambda}
R-r=|y|-|x| \leq 2 {\rm dist}(O,\pi_{\ell_e})=2|\ell_e|\,
\end{equation}
that combined with the fact that ${\rm dist}(O,\pi_{\ell_e}) \leq C \sqrt{\deficit}$ gives
that 
\begin{equation}
\label{viavia}
R-r\leq 2\,C \sqrt{\deficit}\,,
\end{equation}
and by assuming that $\deficit$ is \emph{very} small, we obtain
$$
R-r\leq 2\,C(n) \sqrt{\deficit}\,\sqrt{R-r},
$$
that is
\begin{equation}
\label{eq:Rr delta}
R-r \leq  C(n)\,\deficit
\end{equation}
which gives the optimal estimate for the difference of the radii and we conclude.

\section{Open problems} \label{sect_open}
In this section we state some open problems which are related to quantitative studies of MMP or to the problems which we reviewed in this paper. 

\begin{itemize}
	\item[\emph{Problem 1.}] As mentioned in the Introduction, there are essentially two ways for proving Alexandrov Theorem \ref{thm_Alex}. One is based on the MMP and one involves integral and geometric identities. In the nonlocal setting, the only proof for Theorem \ref{thm_nonlocal} is by using the method of moving planes (see \cite{Cabre} and \cite{CFMN}) and it is not available in literature a proof which uses integral identities, say in the spirit of \cite{MontielRos} or \cite{Ros}. We believe that this is a very interesting and challenging problem.
	
	\item[\emph{Problem 2.}] Once Problem 1 is solved, it would be interesting to understand if it is possible to study the new proof based on integral identities from a quantitative point of view. In the classical setting, the quantitative study of the proof of \cite{Ros} leaded to the quantitative results for the bubbling phenomenon in \cite{CM}. Analogously, in the nonlocal case, a proof based of integral identities might lead to the understanding of the bubbling phenomenon in the nonlocal case by using a deficit different from the one used in \cite{CFMN}.
	
	\item[\emph{Problem 3.}]  Serrin's overdetermined problem in the nonlocal setting was proved in \cite{Jarohs} by using the MMP. However, there are no quantitative results of proximity to a single ball for Serrin's overdetermined problem in the nonlocal case and it is expected that these results may be obtained starting from the proof in \cite{Jarohs} and performing a quantitative analysis of the MMP. 
	
	\item[\emph{Problem 4.}] The same scenario of Problems 1 and 2 occurs also for Serrin's overdetermined problem. More precisely, it is of interest to give a proof of Serrin's overdetermined problem in the nonlocal setting which is not based on the MMP but on integral identities. It is clear that also the quantitative aspects of such a proof would be of strong interest.
	
	\item[\emph{Problem 5.}] Regarding the symmetry problem in Theorem \ref{thm_Rolly} it would be of interest to study the symmetry problem in the nonlocal case as well as its quantitative counterpart. This can be done by using the MMP. Moreover, it would be of interest to provide a proof of Theorem \ref{thm_Rolly} by using an alternative approach, for instance integral identites as it happens for Serrin's overdetermined problem.
	
	\item[\emph{Problem 6.}] The exterior Serrin's overdetermined problem states that if $u$ is the capacitary function of a set $\Omega$ with constant normal derivative on $\partial \Omega$ then $\Omega$ is a ball. This overdetermined problem was solved in \cite{Reichel} by using the MMP (see also \cite{BianchiniCiraolo}, \cite{BianchiniCiraoloSalani} and \cite{SartoriGarofalo} for other approaches). As far as the authors know, a stability result of this theorem is not available in literature; hence a quantitative study of the MMP in this setting would be of interest.
	
	\item[\emph{Problem 7.}] There is an analogue of Alexandrov's Theorem in capillarity theory for embedded constant mean curvature hypersurfaces with boundary contained in a hyperplane with constant contact angle. A proof of this result can be given by using the MMP as well as integral identities (see \cite{MaMi}), and one can start from there to obtain quantitative stability results for this problem.
	
	\item[\emph{Problem 8.}] Probably the most general setting where an Alexandrov's type Theorem has been proved is in warped product manifolds by Brendle \cite{Brendle} (see also \cite{Qui Xia}). The proof is in the spirit of the ones of Montiel and Ros \cite{MontielRos} and Ros \cite{Ros}, and a quantitative investigation of that problem would be of great interest especially in order to characterize the bubbling phenomenon in that context.
	
	\item[\emph{Problem 9.}]  As already mentioned in the Introduction, an Alexandrov's type theorem still holds for immersed hypersurfaces of genus  $0$ (see \cite{Hopf} and \cite{Yau}). The proof does not make use of the MMP and it would be interesting to study a quantitative stability version of this theorem.

	\item[\emph{Problem 10.}] One of the most beautiful applications of the MMP was done in the seminal paper \cite{GNN1}, where it is proved the spherical symmetry of solutions to semilinear equations $\Delta u + f(u)=0$ in $\mathbb{R}^n$. A quantitative study of this problem is a challenging and very interesting problem. It is expected that one can have bubbling and that, in order to have a stability result of proximity to a single ball by using the quantitative MMP one has to introduce some a priori condition on the solution $u$.

\end{itemize}

\bibliographystyle{alpha}

\end{document}